\documentclass[11pt, oneside]{article}

\usepackage{amsmath}
\usepackage{amssymb}
\usepackage{graphics}
\usepackage{graphicx}
\usepackage{theorem}
\usepackage{latexsym}
\usepackage{setspace}

\usepackage{amssymb}

\newcommand{\dd}{\partial\bar{\partial}}
\newcommand{\ip}{\centerdot}

\newtheorem{lemma}{Lemma}
\newtheorem{theorem}{Theorem}

\newtheorem{prop}{Proposition}
\newtheorem{definition}{Definition}

\begin{document}

\title{Spectral measures on toric varieties and the asymptotic expansion of Tian-Yau-Zelditch}
\author{Rosa Sena-Dias}
\maketitle
\begin{abstract}
We extend a recent result of Burns, Guillemin and Uribe on the asymptotics of the spectral measure for the reduction metric on a toric variety to any toric metric on a toric variety. We show how this extended result together with the Tian-Yau-Zelditch asymptotic expansion can be used to deduce Abreu's formula for the scalar curvature of a toric metric on a toric variety in terms of polytope data.
\end{abstract}

\section{Introduction}
Recently, Burns, Guillemin and Uribe have described a procedure to give the asymptotic expansion of the so called spectral measure sequence on a toric manifold. One obtains this sequence by choosing an orthonormal basis for the space of holomorphic sections of large tensor powers of a quantizing line bundle for the toric manifold and adding the square of the norms of these elements of the basis. More precisely, suppose we are given a toric symplectic manifold $(X,\omega)$. Let $L\rightarrow X$ be a line bundle whose first Chern class is $[\omega]$. There is a natural metric on $X$ called the reduction metric which is invariant by the torus action on $X$ (see \cite{g1}). This metric allows us to define a Hermitian metric on the bundle $L$. We simply ask that this Hermitian metric $h$ satisfies the following: $i\dd \log h$ is the K\"ahler metric on $X$. The space of holomorphic sections of $L^N$ which we denote by $H^0(X,L^N)$ inherits a torus action and splits according to the action's weights. Since we have a Hermitian metric on $L$, we can choose an orthonormal basis $\{s_m\}$ of $H^0(X,L^N)$ which is compatible with the splitting of $H^0(X,L^N)$. The sequence of spectral measures is
\begin{displaymath}
 \mu_N=\sum_i |s_{m}|^2\nu
\end{displaymath}
where $\nu$ is the Liouville measure. One of the results in \cite{bgu} describes the asymptotic behavior of the sequence $\mu_N$ in $N$. 

It is well known (see \cite{g1}) that toric varieties have many toric metrics (i.e. metrics invariant through the torus action) compatible with a given symplectic structure. Each toric metric gives rise to a different Hermitian metric on the bundle $L$ hence to a different orthonormal basis from which it is possible to construct a spectral measure depending on the initial toric metric. It is then natural to ask if the results in \cite{bgu} extend to such metrics. The first purpose of this note is to show that this is true. For example, we prove 
\begin{theorem}\label{main}
 Let $X$ be a toric variety with moment polytope $\triangle$ and consider any toric metric on $X$. Let $\psi \in \mathcal{C}^\infty(\triangle)$, and $\mu_N$ be the sequence of spectral measures on $X$ for the chosen toric metric. Then
 \begin{displaymath}
 \int\psi \mu_N=\sum_{i=0}^NP_i(\psi)N^{n-i}
\end{displaymath}
where $\psi$ can be seen as a function on $X$ via the moment map and the $P_i(\psi)$'s are integrals of differential operators acting on $\psi$.
\end{theorem}
The pointwise asymptotic behavior of the sequence of functions
\begin{displaymath}
 \sum |s_m|^2
\end{displaymath}
has been extensively studied (see \cite{catlin}, \cite{lu}, \cite{t} and \cite{z}). It is known as the Tian-Yau-Zelditch asymptotic expansion. Now when one applies the measures $\mu_N$ to functions with compact support on the polytope, one should be able to recover the pointwise asymptotic behavior of the function $\sum |s_m|^2$ (at least in ``most'' points). In this spirit, we use Theorem \ref{main} in a precise version (or rather the method that is used to prove this theorem and which appears in \cite{bgu})) to write down the term in $N^{n-1}$ in the expansion. By comparing this term with the corresponding term of the pointwise Tian-Yau-Zelditch asymptotic expansion (which was obtained by Lu in \cite{lu}) we recover a result of Abreu's (see \cite{m1}) which gives a formula for the scalar curvature of any toric metric on a toric variety in terms of polytope data. That is, we give a different proof of the following theorem of Abreu
\begin{theorem}
 Let $X$ be a toric variety with moment polytope $\triangle$ and consider any toric metric on $X$ with symplectic potential $g$. The scalar curvature of the metric is given by
 \begin{displaymath}
 -\frac{1}{2}\frac{\partial^2G^{ij}}{\partial y_i \partial y_j}
\end{displaymath}
where $y$ is the moment map coordinate for $g$ on $\triangle$ and $G$ is the Hessian of $g$ in the $y$ coordinates.

\end{theorem}

A brief outline of this paper is the following: In the first section we give a very brief review of some old results on the K\"ahler structure of toric varieties. The results described in that section first appeared in \cite{g1} and \cite{m1}. The second section deals with the function $\varphi$ which encodes the information coming from the norms of a special basis of $H^0(X,L^N)$. In this section, we basically prove that the function $\varphi$ satisfies properties similar to the ones appearing in \cite{bgu} which ensure that the sections $s_m$ peak at some fiber over a point in $\triangle$. In the fourth section, we use these properties to write down the asymptotic behavior of $|s_m|$ and of $\mu_N$. Finally, in the fifth section we compare these results with the well known Tian-Yau-Zelditch asymptotic expansion to get a formula for the scalar curvature of the toric metric in terms of polytope data only. This formula was first obtained in \cite{m1} using a different approach.

{\bf Acknowledgments} Most of this work was carried out while I was at IST. I would like to take the opportunity to thank Miguel Abreu for his support, for many interesting discussions on K\"ahler metrics on toric manifolds and for having introduced me to the subject some years ago. I am also very grateful to Victor Guillemin for his enthusiasm for this work and for many valuable suggestions on a preliminary version of this paper.

\section{Background}
For the sake of completeness we give some background on K\"ahler toric varieties (see \cite{g1}, \cite{m1} and \cite{m2} for more details and proofs).
\subsection{K\"ahler toric metrics}

\begin{definition}
 A K\"ahler toric manifold $X^{2n}$ is a closed connected K\"ahler manifold $(X,\omega,J)$ with a $T^n$ Hamiltonian action which is also holomorphic.
\end{definition}
Such an action admits a moment map $\phi:X \rightarrow \mathbb{R}^n$, where we have identified $\mathbb{R}^n$ with its dual. This moment map depends on $\omega$ but its image does not. It is a convex polytope in $\mathbb{R}^n$ of Delzant type.
\begin{definition}
 A convex polytope $\triangle$ in $\mathbb{R}^n$ is of Delzant type if
 \begin{enumerate}
 \item There are $n$ edges meeting at each vertex,
 \item It is possible to choose a set of primitive exterior normals to facets in $\mathbb{Z}^n$,
 \item The set of outward normals corresponding to the facets meeting at one vertex forms a basis for $\mathbb{R}^n$, for all vertexes.
\end{enumerate}
\end{definition}
Given a polytope in $\mathbb{R}^n$, Delzant has given a canonical way to associate to it a K\"ahler toric manifold which we write $(X,\omega_0,J)$. The results in \cite{bgu} concern this K\"ahler manifold. The symplectic form $\omega_0$ and the metric associated with $\omega_0$ and $J$ are called reduced.

Although all K\"ahler toric manifolds are symplectomorphic to $(X,\omega_0)$ they are by no means K\"ahler isomorphic to $(X,\omega_0,J)$. That is, the symplectic isomorphism between $(X,\omega)$ and $(X,\omega_0)$ can, in general, not be taken to be holomorphic. In a dual way, we can also say that all K\"ahler toric manifolds are biholomorphic to $(X,J)$ but the biholomorphism can not be taken to be symplectic. Therefore we can think of K\"ahler toric manifolds as triples $(X,\omega,J)$ where $J$ is canonical and comes from the Delzant construction but $\omega$ is, in general, different from $\omega_0$.

\subsection{Complex and symplectic coordinates}
Toric manifolds admit two natural sets of coordinates (see \cite{m2}):
\begin{enumerate}
\item Complex coordinates on an open dense subset (the open orbit of the torus action). Let us call these coordinates $z=u+iv$. In these coordinates the torus acts by translation in the $v$ variables. 

\item There are also symplectic coordinates, on the same dense open set of $X$, that make the symplectic form standard. Consider the moment map on $X$ from the torus action. As is well known, its image lies in a Delzant type polytope in $\mathbb{R}^n$ which we call $\triangle$. Let $y$ be the moment map coordinates on this polytope, then $(y,v)$ are symplectic coordinates in $X$, that is
\begin{displaymath}
 \omega=\sum dy_i\wedge dv_i.
\end{displaymath}
In these coordinates, the complex structure takes a non-standard form
\begin{displaymath}
 \begin{pmatrix}
  0& hess{g}\\- hess{g}&0
 \end{pmatrix}.
\end{displaymath}
for some function $g$ which is called the symplectic potential. Note that, since the complex structure is torus invariant, $g$ is really only a function of the $y$ coordinates.
\end{enumerate}

There is a relation between complex and symplectic coordinates, namely the Legendre transform
\begin{equation}\label{Legendre}
u=g_y.
\end{equation}
Also in the open dense subset, the K\"ahler form admits a K\"ahler potential. That is, there exits $f$, a function of $z$, such that the K\"ahler form is $-2i\dd f$. The invariance of the K\"ahler metric with respect to the torus action implies that $f$ is actually only a function of $u$. The functions $f$ and $g$ are related to each other by the following relations (see \cite{g1})
\begin{equation}\label{fg}
 f(u)+g(y)=y \ip u, \qquad y=f_u, \qquad u=g_y.
\end{equation}
\subsection{General symplectic potentials}\label{generalpotentials}
The polytope $\triangle$ can be described by a set of inequalities
\begin{displaymath}
 \triangle=\{y\in \mathbb{R}^n: y\ip u_i-c_i \leq 0, i=1,\cdots, d\}
\end{displaymath}
where $c_i$ is an integer and $u_i \in \mathbb{Z}^n$ is a primitive outward normal to the $i$-th facet of $\triangle$. Set $l_i(y)=c_i-y\ip u_i$. The reduction metric corresponding to $\omega_0$ has symplectic potential
\begin{displaymath}
 g_0=\sum_{i=1}^d l_i\log l_i-l_i.
\end{displaymath}
The symplectic potential of a general metric is given by $g=g_0+g_r$ where $g_r$ is smooth in a neighborhood of $\triangle$ (see \cite{m2}). In particular the behavior of the potential on the facets of the polytope is that of $g_0$. 

Not all functions of the form $g=g_0+g_r$, where $g_r$ is smooth, are symplectic potentials though (see Theorem 2.8 in \cite{m2}). In particular, $g$ needs to be such that $hess(g)$ is positive definite on the interior of $\triangle$ and this Hessian needs to have a certain behavior on the facets. For example, as one reaches a point in the interior of an $(n-1)$-dimensional facet, the inverse matrix of the Hessian converges but acquires a kernel which is generated by the $u_i$ corresponding to that facet. As one approaches an $(n-2)$-dimensional facet, the inverse of the Hessian still converges but this time acquires a $2$-dimensional kernel etc... 

\subsection{Holomorphic sections on toric manifolds}
The $T^n$ action on $X$ is holomorphic and therefore it induces an action on $H^0(X,L^N)$, the space of holomorphic sections on $L^N$. This vector space must split according to the weights of the action. There is a torus invariant global section which we denote by $\mathbf{1}$. It is not hard to see that any holomorphic section of $L^N$ can be written as a linear combination of the sections $e^{m \ip   z}\mathbf{1}^N$ i.e.
\begin{displaymath}
 H^0(X,L^N)=span\{e^{m \ip   z}\mathbf{1}^k, \quad m\in \mathbb{Z}^n \cap N\triangle,  \}.
\end{displaymath}
We set $\mathbb{Z}^n\cap N\triangle=[N\triangle]$. This basis decomposes  $H^0(X,L^N)$ into one dimensional weight spaces for the torus action. Namely,
\begin{displaymath}
 e^{i\theta}e^{m \ip   z}=e^{i m\ip\theta}e^{m \ip   z},
\end{displaymath}
 where $\theta$ is in $\mathbb{R}^n$ and we write $e^{m \ip   z}$ for $e^{m \ip   z}\mathbf{1}^N$. The set $\{e^{m \ip   z}, \quad m\in \mathbb{Z}^n \cap N\triangle,  \}$, forms an orthogonal basis of $H^0(X,L^N)$. This is simply because
\begin{displaymath}
 \int_{T^n}e^{im\ip v}dv=0,
\end{displaymath}
unless m=0.

\section{The function $\varphi$}
Our setting is the same as the setting in \cite{bgu}. Let $X$ be a K\"ahler toric manifold of complex dimension $n$ such that the symplectic form on $X$, $\omega$ has integral cohomology class so that it is the curvature of some connection on a line bundle $L\rightarrow X$. Now consider any toric metric on $X$, i.e. any metric which is invariant by the torus action on $X$. This metric comes from a Hermitian metric on $L$, h. In fact, the K\"ahler metric on $X$ is given by $i\dd \log h$.
Pick an orthonormal basis for $H^0(X,L^N)$ which is an eigenbasis for the torus action, say $\{s_{m}\}$. As in \cite{bgu}, we are interested in the asymptotic behavior of the spectral measure
\begin{displaymath}
 \mu_N=\sum_i |s_{m}|^2\nu
\end{displaymath}
where $\nu$ is the Liouville measure. In \cite{bgu}, Burns, Guillemin and Uribe consider the case where the metric on $X$ is the so called reduced metric. Here we are concerned with the general case.
To this end, we look at the norms of the sections $e^{m \ip   z}$ with respect to the Hermitian metric associated to $\omega$ and thus define the function $\varphi$ which encodes the information from all of these norms. 
\subsection{Definition}  
Suppose we have a fixed toric metric on $X$. This metric allows us to define a Hermitian metric on the bundle $L$. Simply set
\begin{displaymath}
 \omega=i\dd \log h
\end{displaymath}
where $h$ is the norm of the torus invariant section we have called $\mathbf{1}$. The function $h$ is torus invariant. We note here that the norm of any of the $e^{m \ip   z}$  is also torus invariant. This is because 
\begin{displaymath}
 |e^{m \ip   z}|_h^2=e^{2m \ip   u}h,
\end{displaymath} 
which does not depend on the $v$ coordinate. We define
\begin{displaymath}
 \varphi\left(\frac{m}{N},y\right)=\frac{1}{2N}\log |e^{m\ip z}|_{h}^2\circ \phi^{-1}(y),
\end{displaymath}
where $m\in \mathbb{Z}^n \cap N\triangle$ and $\phi$ is the moment map for the torus action with respect to $\omega$. Even though $\phi^{-1}$ is not a well defined function, $\varphi$ is well defined. To be more precise, suppose that two points in $X$ have the same image via the moment map in the polytope. Then, they are in the same torus orbit and therefore, since $|e^{m\ip z}|_{h}^2$ is torus invariant, the above quantity is well defined. In \cite{bgu}, Burns, Guillemin and Uribe consider the case where the symplectic form is the reduction symplectic form. Then we have
\begin{displaymath}
 \varphi_0\left(\frac{m}{N},y\right)=\frac{1}{2N}\log |e^{m\ip z}|_{h_0}^2\circ \phi_0^{-1}(y),
\end{displaymath}
where $h_0$ is the Hermitian metric corresponding to the reduction symplectic form and $\phi_0$ is the moment map of the torus action associated with the reduction symplectic form.
\subsection{$\varphi$ and $\varphi_0$}
We can write down a relation between the functions $\varphi$ and $\varphi_0$. For that we need to consider the map
\begin{displaymath}
 \alpha(y)=\phi_0\circ \phi^{-1}(y).
\end{displaymath}
 This is well defined because if two points have the same image via $\phi$ then they lie in the same torus orbit above that image and therefore they have the same image via $\phi_0$. Note that in the interior of $\triangle$ there an explicit expression for $\alpha$, namely
\begin{displaymath}
 \alpha(y)=f_{0u}\circ f_u^{-1}.
\end{displaymath} 
For later use we prove the following simple lemma
\begin{lemma}\label{alpha}
The function $\alpha$ is the identity on the boundary of $\triangle$.
\end{lemma}
{\bf Proof} The proof is simple. The point is that the function $f_{0u}\circ f_u^{-1}$ extends smoothly to the boundary of $\triangle$ as the identity. The function $f_{0u}\circ f_u^{-1}$ is simply ${g_{0y}}^{-1}\circ g_y$, and the behavior of $g_y$ is the same as that of $g_{0y}$ at the boundary of $\triangle$. To be more precise let $x$ be in the boundary of $\triangle$, $I$ be set of faces to which $x$ does not belong and $I^c$ its complement in the set $\{1, \cdots,d\}$ i.e. the set of indices of the facets to which $x$ does belong.
\begin{displaymath}
 g_y(y)=-\sum_{i=1}^d \log l_i(y) u_i+g_{ry}\simeq  \sum_{i\in I^c} \log l_i(y)u_i,
\end{displaymath}
near the interior of the $i$th face of $\triangle$. Since the same is true for $g_{0y}$, the result follows.$\square$

The forms $\omega_0$ and $\omega$, are cohomologous therefore there is a globally defined function on $X$, say $\rho$, such that
\begin{displaymath}
 \omega=\omega_0+2\dd \rho.
\end{displaymath}
Seen as a function of $z$ in the open dense orbit, $\rho$ only depends on $u$ because it must be invariant by the torus action.
\begin{lemma}
The functions $\varphi$ and $\varphi_0$ are related via
\begin{displaymath}
 \varphi(x,y)=\varphi_0(x,\alpha(y))+\rho(\alpha(y))
\end{displaymath}
\end{lemma}
{\bf Proof} This is straightforward.  We must have
\begin{displaymath}
 h=e^{2\rho}h_0,
\end{displaymath}
up to a constant, hence 
\begin{displaymath}
 |e^{m\ip z}|_{h}^2=|e^{m\ip z}|_{h_0}^2e^{2N\rho},
\end{displaymath}
and
\begin{displaymath}
 \varphi(x,y)=\frac{1}{2N}\log |e^{m\ip z}|_{h}^2\circ \phi_0^{-1}(\phi_0\circ \phi^{-1}(y))
\end{displaymath}
where $x=m/N$. That is 
\begin{displaymath}
 \varphi(x,y)=\frac{1}{2N}\log\left(|e^{m\ip z}|_{h_0}^2e^{2N\rho}\right)\circ \phi_0^{-1}(\alpha(y))
\end{displaymath}
and the result follows.$\square$
\subsection{$\varphi$ and $g$}
Let $g$ denote, as before, the symplectic potential associated to the symplectic form $\omega$. It is possible to write down $\varphi$ in terms of $g$ alone.
\begin{lemma}
Let $g$ be the symplectic potential for $\omega$, we have 
 \begin{displaymath}
 \varphi(x,y)=g(y)+(x-y)\ip g_y(y).
\end{displaymath}
\end{lemma}
{\bf Proof} We start by determining an expression for the Hermitian metric $h$. Since both $-2f$ and $\log h$ are potentials for the K\"ahler metric on $X$ we must have
\begin{displaymath}
 h=e^{-2f}.
\end{displaymath}
Therefore
\begin{displaymath}
 |e^{m\ip z}|_{h}^2=e^{2m\ip u}e^{-2Nf}.
\end{displaymath}
Replacing $f$ by the expression given in equation (\ref{fg}) and using equation (\ref{Legendre}) we have
\begin{displaymath}
 |e^{m\ip z}|_{h}^2(y)=e^{2N(g(y)+(\frac{m}{N}-y)\ip g_y(y))},
\end{displaymath}
where $y=\phi(z)$. We have
\begin{displaymath}
 \varphi\left(\frac{m}{N},y\right)=\frac{1}{2N}\log |e^{m\ip z}|_{h}^2\circ \phi^{-1}(y)
\end{displaymath}
and the result follows.$\square$

Let us check that this fits in well with the expression in \cite{bgu} for $\varphi_0$. Since
\begin{displaymath}
 g_0=\sum_{i=1}^d l_i\log l_i-l_i,
\end{displaymath}
we have
\begin{displaymath}
 g_{0y}=-\sum_{i=1}^d  u_i \log l_i,
\end{displaymath}
and therefore
\begin{displaymath}
 g_0+(x-y)\ip g_{0y}=\sum l_i(x)\log l_i(y)-l_i(y),
\end{displaymath}
because $x-y \ip u_i=l_i(x)-l_i(y)$. The above expression coincides with the expression appearing in \cite{bgu} for $\varphi_0$. 
 
Later we are going to be interested in an orthonormal basis for $H^0(X,L^N)$. This is simply the set $\{s_m\}$, where
\begin{displaymath}
 s_m=\frac {e^{m\ip z}}{||e^{m\ip z}||}
\end{displaymath}
where $||.||$ refers to the $L^2$ norm and we can write
\begin{displaymath}
 |s_m|_{h}^2 (y)=\frac {e^{2N\varphi(\frac{m}{N},y)}} {\int_{\triangle}e^{2N\varphi(\frac{m}{N},y)} dy}, \quad m\in [N\triangle].
\end{displaymath}
\subsection{Two lemmas on $\varphi$}
The point here is that independently of $g$, the function $\varphi$ satisfies two lemmas which appear in \cite{bgu} for the case $g=g_0$.
\begin{lemma}\label{maximo}
 Let $x$ be a point in the interior of $\triangle$ and $g$ be a symplectic potential on $\triangle$. Then, the function $\varphi$ regarded as a function of $y=f_u(u)$ has a unique critical point at $x=y$ and this unique critical point is the unique global maximum of the function $\varphi$ on $\triangle$.
\end{lemma}
{\bf proof } We first note that as $y$ tends to $\partial \triangle$, $\alpha(y)$ also tends to $\partial \triangle$. Consider the formula 
\begin{displaymath}
 \varphi(x,y)=\varphi_0(x,\alpha(y))+\rho(\alpha(y)).
\end{displaymath}
We know that $\varphi_0(x,y)$ tends to $-\infty$ as $y$ tends to $\partial \triangle$ because  
\begin{displaymath}
 \varphi_0(x,y)=\sum l_i(x)\log l_i(y)-l_i(y).
\end{displaymath}
As for $\rho$, since it a globally defined function on $X$ it must have a finite limit as $y$ tends to $\partial \triangle$. We conclude that $\varphi(x,y)$ tends to $-\infty$ as well on $\partial \triangle$. On the other hand it is bounded from above on $\triangle$ since the $l_i$ and $\rho$ are. Therefore it has a maximum on the interior on $\triangle$. This maximum is a critical point of $\varphi$. Using
\begin{displaymath}
 \varphi(x,y)=g(y)+(x-y)\ip g_y(y), 
\end{displaymath} 
we see that
\begin{equation}\label{derivadaphi}
 \frac{\partial\varphi}{\partial y}=hess(g) (x-y).
\end{equation}
Now from the properties of $g$ mentioned in section \ref{generalpotentials} we know that $hess(g)$ is positive definite on the interior of $\triangle$ and the result follows.$\square$

A similar result can be proved when $x$ is in the boundary of $\triangle$ namely
\begin{lemma}\label{maximo_fronteira}
 Let $x$ be a point in the interior of a facet $F$ of $\triangle$ and $g$ be a symplectic potential on $\triangle$. Then, the restriction to $F$ of the function $\varphi$ regarded as a function of $y$ has a unique critical point at $x=y$ and this unique critical point is the unique global maximum of the restriction to $F$ of the function $\varphi$ on $\triangle$. Moreover, the derivatives of $\varphi$ in the directions normal to $F$ are not zero at this maximum. 
\end{lemma}
{\bf proof} Let $I=\{i\in \{1,\cdots,d\}: x\notin l_i^{-1}(0)\}$ where $d$ is the total number of facets in $\triangle$. That is, $I$ is the set of indices of the facets to which $x$ does not belong. We start by showing that $\varphi(x,\cdot)$ actually extends to $F$. This follows from the formula
\begin{displaymath}
 \varphi(x,y)=\sum_{i\in I} l_i(x)\log l_i(\alpha(y))-\sum_{i=1}^d l_i(\alpha(y))+\rho(\alpha(y))
\end{displaymath}
which in turn is a consequence of the expression
\begin{displaymath}
 \varphi_0(x,y)=\sum_{i\in I} l_i(x)\log l_i(y)-\sum_{i=1}^d l_i(y).
\end{displaymath}
from \cite{bgu}. 
Again as $y$ tends to $\partial F$, so does $\alpha(y)$ and it follows from the expression above that $\varphi_0(x,\cdot)$ tends to $-\infty$. Therefore $\varphi(x,\cdot)$ is $-\infty$ on the boundary of $F$. It is also bounded from above on this facet so there must be a maximum on the interior of $F$ and this maximum must be a critical point of the restriction of $\varphi$ to $F$ as a function of $y$. Consider the expression 
\begin{equation}\label{x-y}
 x-y=(hess_y(g))^{-1}\frac{\partial\varphi}{\partial y}
\end{equation}
which holds true on the interior of $\triangle$. We know from the properties of $g$ described in section \ref{generalpotentials} and discussed in $\cite{m2}$ that $(hess_y(g))^{-1}$ extends to $F$ with a kernel generated by $\{u_i, i\in I^c\}$. As for $\frac{\partial\varphi}{\partial y}$ it is given by
\begin{displaymath}
 \frac{\partial\varphi}{\partial y}(y)= D \alpha(y)\left(\frac{\partial\varphi_0}{\partial y}(x,\alpha(y))+\frac{\partial\rho}{\partial y}(\alpha(y))\right)
\end{displaymath}
and 
\begin{displaymath}
\frac{\partial\varphi_0}{\partial y}=-\sum_{i\in I} \frac{l_i(x)}{l_i(y)}u_i+\sum_{i=1}^d u_i
\end{displaymath}
which clearly extends to the interior of $F$. Hence we conclude that $\frac{\partial\varphi}{\partial y}$ itself extends to the interior of $F$. So equation (\ref{x-y}) holds even for $y$ in $F$. Suppose that the point $y\in F$ is critical for the restriction of $\varphi$ to $F$. This means that 
\begin{displaymath}
 \frac{\partial\varphi}{\partial y} \in T^{\perp}F=span \{u_i, i\in I^c\}
\end{displaymath}
hence
\begin{displaymath}
 (hess_y(g))^{-1}\frac{\partial\varphi}{\partial y}=0
\end{displaymath}
which implies that $x=y$. Next we would like to see that $\frac{\partial\varphi}{\partial y}$ cannot be zero at $x=y$. Define
\begin{displaymath}
 g_{0I^c}=\sum_{i \in I^c}l_i\log l_i -l_i.
\end{displaymath}
We have $g=g_{0I^c}+g_0-g_{0I^c}+g_r$. The functions $g_r$ and $g_0-g_{0I^c}$ extend smoothly to the interior of $F$ so we must check that
\begin{displaymath}
 hess_y(g_{0I^c})(x-y)
\end{displaymath}
tends to something which is not zero as $y$ tends to $x$. We have
\begin{displaymath}
 hess(g_{0I^C})(y)=\sum_{i \in I^c}\frac{u_i^t u_i}{l_i(y)}
\end{displaymath}
hence
\begin{displaymath}
 hess(g_{0I^C})(y)(x-y)=\sum_{i \in I^c}\frac{ u_i\ip (x-y)}{l_i(y)}u_i.
\end{displaymath}
Now use the fact that $ u_i\ip (x-y)=l_i(y)-l_i(x)=l_i(y)$ so the expression above coincides with
\begin{displaymath}
\sum_{i \in I^c} u_i.
\end{displaymath}
which is non zero because of the Delzant condition.$\square$
\section{The results}
As in \cite{bgu} we are interested in the spectral measures of the manifold $X$. These are defined by
\begin{displaymath}
 \mu_N=\sum |s_m|_h^2\nu
\end{displaymath}
where $\nu$ is the Liouville measure. Let $\psi$ be a smooth function on $\triangle$. Then the goal is to write an asymptotic formula for 
\begin{displaymath}
 \int_\triangle \psi \mu_N,
\end{displaymath}
in $N$, as $N$ tends to $\infty$. We have 
\begin{displaymath}
 \int_\triangle \psi \mu_N=\sum_{m \in [N\triangle]} \psi^\sharp \left( \frac{m}{N} \right)
\end{displaymath}
where
\begin{displaymath}
 \psi^\sharp(x)=\int_{\triangle}\frac{\psi(y)e^{2N\varphi(x,y)}}{\int_{\triangle}{e^{2N\varphi(x,y)}}} dy.
\end{displaymath}
The results in \cite{bgu} use two main ingredients:
\begin{itemize}
 \item The first in an asymptotic formula for sums of the form
 \begin{displaymath}
 \sum_{m\in [N\triangle]} \psi\left(\frac{m}{N}\right)
\end{displaymath}
for any continuous function $\psi$ on $\triangle$ (see \cite{gs}). This argument clearly does not depend on the symplectic potential $g$.
\item The second ingredient is an application of the Euler-McLaurin formula to the function $\psi$ as it appears in the integral 
\begin{displaymath}
 \int_{\triangle}\frac{\psi(y)e^{2N\varphi(x,y)}}{\int_{\triangle}e^{2N\varphi(x,y)} dy}dy
\end{displaymath}
around the point $x$ which is the point where $\varphi$ attains its maximum. Again this argument works for general $g$ since the lemmas 4.1 and 4.2 of \cite{bgu} carry over to this case. Their generalizations to this setting are lemmas \ref{maximo} and \ref{maximo_fronteira} from the previous section. We will carry out this method explicitly in the next section for the case when $\psi$ has compact support in $\triangle$.
\end{itemize}
We can summarize the results obtained by applying this method in the following generalization of the main theorem in \cite{bgu}:
\begin{theorem}
 Let $X$ be a toric variety with moment polytope $\triangle$ endowed with a toric metric. Let $\psi \in \mathcal{C}^\infty(\triangle)$, and $\mu_N$ be the sequence of spectral measures on $X$ for the chosen toric metric. Then
 \begin{displaymath}
 \int\psi \mu_N=\sum_{i=0}^NP_i(\psi)N^{n-i}
\end{displaymath}
where $\psi$ can be seen as a function on $X$ via the moment map and the $P_i(\psi)$'s are integrals of differential operators acting on $\psi$.
\end{theorem} 
The other asymptotic results appearing in \cite {bgu} hold true in this new setting provided one is careful to note that the coordinates $y$ are now given by the new moment map $\phi$ corresponding to $\omega$. We have for example:
\begin{theorem}
 Let $x$ be in $\triangle$ and $m=Nx$. The sequence of sections $s_m$ converges to a delta function on the fiber $\phi^{-1}(x)$.
\end{theorem}
{\bf proof } This is exactly as in \cite{bgu}. We note that
\begin{displaymath}
 \int_{\triangle}e^{N\varphi(x,y)}dy\simeq \left(\frac{2\pi}{N}\right)^{n/2}h(x)^{-1/2}e^{N\varphi(x,x)},
\end{displaymath}
where $h(x)$ is the determinant of $hess(g)$. Therefore, 
\begin{displaymath}
 |s_m|^2\simeq \left(\frac{N}{2\pi}\right)^{n/2}h(x)^{1/2}e^{N(\varphi(x,y)-\varphi(x,x))}.
\end{displaymath}
Since for $y$ not equal to $x$ we have $\varphi(x,y)-\varphi(x,x)<0$, the above converges to zero except if $x=y$, that in on the $x$ fiber of $\phi$.$\square$

\section{Abreu's scalar curvature formula}
\subsection{An explicit calculation using \cite{bgu} approximation method}
In the previous section we roughly described the method first presented in \cite{bgu} to obtain the asymptotic behavior of the integral
\begin{displaymath}
 \int_{\triangle}\psi\mu_N.
\end{displaymath}
In the case where $\psi$ has compact support in $\triangle$, the method simplifies considerably. We are going to write down explicitly the first two terms in the expansion and see how, from the second term, we can recover Abreu's formula (see \cite{m1}) for the scalar curvature of a toric manifold.

Let $\psi$ be in $\mathcal{C}_0^\infty(\triangle)$. First write 
\begin{displaymath}
 \int_{\triangle}\psi\mu_N=\sum_{m\in [N\triangle]} \psi^\sharp(\frac{m}{N}).
\end{displaymath}
From \cite{gs} it is known that
\begin{equation}\label{somapolytopo}
 \frac{1}{N^n}\sum_{m\in [N\triangle]} \psi^\sharp(\frac{m}{N})\sim \tau(\frac{1}{N}\frac{\partial}{\partial h})\int_{\triangle_h}\psi^\sharp(h=0),
\end{equation}
where $h\in \mathbb{R}^d$. By $\triangle_h$, we mean the dilated polygon 
\begin{displaymath}
 \triangle_h=\{y\in \mathbb{R}^n: y\ip u_i-c_i \leq h_i, i=1,\cdots, d\}.
\end{displaymath}
The function $\tau$ is defined by
\begin{displaymath}
 \tau(s):=\frac{s}{1-e^s}=1+\frac{s}{2}+O(s^2)
\end{displaymath}
and
\begin{displaymath}
 \tau(\frac{1}{N}\frac{\partial}{\partial h})=\tau(\frac{1}{N}\frac{\partial}{\partial h_1})\cdots\tau(\frac{1}{N}\frac{\partial}{\partial h_d}),
\end{displaymath}
so that 
\begin{displaymath}
 \tau(\frac{1}{N}\frac{\partial}{\partial h})=1+\frac{1}{N}\left(\frac{\partial}{\partial h_1}+\cdots+\frac{\partial}{\partial h_d} \right)+O\left(\frac{1}{N^2}\right).
\end{displaymath}
We now move on to write the asymptotics for 
\begin{equation}\label{sharp}
 \psi^\sharp(x)=\int_{\triangle}\frac{\psi(y)e^{2N\varphi(x,y)}}{\int_{\triangle}e^{2N\varphi(x,y)} dy}dy.
\end{equation}
\begin{prop}
The first terms in the asymptotic expansion for $\psi^\sharp$ are given by
 \begin{displaymath}
 \psi^\sharp(x)=\psi(x)+\frac{1}{2N}\left( \frac{1}{2}\frac{\partial^2 \psi}{\partial y_a\partial y_b}G^{ab}+\frac{\partial\psi}{\partial y_a}\frac{\partial G^{ab}}{\partial y_b}\right)+O\left(\frac{1}{N^2}\right).
\end{displaymath}
\end{prop}
{\bf proof} We can write the Taylor expansion for $\psi$ around $x$
\begin{displaymath}
 \psi(y)=\psi(x)+\frac{\partial\psi}{\partial y_a}(\xi(x,y))(y_a-x_a)
\end{displaymath}
where $\xi$ is a smooth function satisfying
\begin{displaymath}
 \xi(x,y) \in \bar{xy},
\end{displaymath}
\begin{displaymath}
 \xi(x,x)=x,
\end{displaymath}
and 
\begin{displaymath}
 \xi(x,y)=\xi(y,x).
\end{displaymath}
Here, $\bar{xy}$ denotes the set $\{tx+(1-t)y, t\in [0,1]\}$. These properties imply that 
\begin{displaymath}
 \frac{\partial \xi_a}{\partial y_b}(x,y)=\frac{\delta_{ab}}{2}.
\end{displaymath}
From equation (\ref{derivadaphi}) we can write
\begin{displaymath}
 y-x=-(hess(g))^{-1}\frac{\partial \varphi}{\partial y}
\end{displaymath}
or, writing $(hess(g))^{-1}=(G^{ab})$
\begin{displaymath}
 y_a-x_a=-G^{ab}\frac{\partial \varphi}{\partial y_b},
\end{displaymath}
so that in the Taylor expansion for $\psi$ we can write 
\begin{displaymath}
 \psi(y)=\psi(x)-\frac{\partial\psi}{\partial y_a}(\xi(x,y))G^{ab}\frac{\partial \varphi}{\partial y_b}.
\end{displaymath}
Now 
\begin{displaymath}
 \frac{\partial \varphi}{\partial y_b}e^{2N\varphi(x,y)}=\frac{1}{2N}\frac{\partial e^{2N\varphi(x,y)}}{\partial y_b}.
\end{displaymath}
We can replace for $\psi$ in equation  (\ref{sharp}) and write 
\begin{displaymath}
 \psi^\sharp(x)=\psi(x)-\int_{\triangle}\frac{1}{\int_{\triangle}e^{2N\varphi(x,y)} dy}\frac{\partial\psi}{\partial y_a}(\xi(x,y))G^{ab}\frac{1}{2N}\frac{\partial e^{2N\varphi(x,y)}}{\partial y_b}dy.
\end{displaymath}
Next we integrate by parts. We do not pick any boundary terms. Note that even for $y$ in the boundary of $\triangle$, $\xi(x,y)$ is not necessarily in the boundary so that $\psi(\xi(x,y))$ is not necessarily $0$. But the term $G^{ab}$ does vanish on the $b$th facet of $\triangle$. The easiest way to see this is to choose coordinates so as to standardize the $bth$ facet to have normal $u_b=e_b$. Then the boundary behavior of $G^{ab}$ implies that $G^{ab}u_b=0$ on the $b$th facet and therefore $G^{ab}=0$ on the $b$th facet. Note also that for fixed $b$ we need only to integrate by parts in the variable $b$. 
\begin{lemma}
In integrating
\begin{displaymath}
  \int_{\triangle}\frac{\partial\psi}{\partial y_a}(\xi(x,y))G^{ab}\frac{\partial e^{2N\varphi(x,y)}}{\partial y_b}dy.
\end{displaymath}
by parts, we do not pick any boundary terms.
\end{lemma}
{\bf proof} Fix $b$. We choose coordinates $y_1,\cdots, y_n$ centered at one of the vertexes in $F_b$ so that 
\begin{itemize}
 \item In these coordinates, $\triangle$ is a subset of the positive octant of $\mathbb{R}^n$,
 \item For all $i=1,\cdots, n$, each facet $F_i=l_i^{-1}(0)$ is contained in the set $\{(y_1,\cdots,y_n):y_i=0\}$.
\end{itemize}
Let $Q$ denote a rectangle in  $\mathbb{R}^n$ of the form $[0,\alpha_1 ]\times \cdots \times [0,\alpha_n ]$ containing $\triangle$. One can find such a $Q$ as long as the $\alpha_i$'s are big enough. Extend $\psi$ by zero to all of $Q$. This extension is still smooth because the support of $\psi$ is contained in the interior of $\triangle$. We will also need to extend $G^{ab}$ to $Q$ in a smooth way so that it is zero when one of the $y_b$'s is $\alpha_b$. Now
\begin{equation}\label{integral}
 \int_{\triangle}\frac{\partial\psi}{\partial y_a}(\xi(x,y))G^{ab}\frac{\partial e^{2N\varphi(x,y)}}{\partial y_b}dy
\end{equation}
is the same as 
\begin{displaymath}
 \int_{Q} \frac{\partial\psi}{\partial y_a}(\xi(x,y))G^{ab}\frac{\partial e^{2N\varphi(x,y)}}{\partial y_b}dy
\end{displaymath}
which we write as 
\begin{displaymath}
 \int_{\widehat{Q}_b}\int_0^{\alpha_b}\frac{\partial\psi}{\partial y_a}(\xi(x,y))G^{ab}\frac{\partial e^{2N\varphi(x,y)}}{\partial y_b}dy_bd\hat{y}
\end{displaymath}
where 
\begin{displaymath}
 \widehat{Q}_b=[0,\alpha_1 ]\times \cdots\times \widehat{[0,\alpha_b]} \times\cdots \times [0,\alpha_n ],
\end{displaymath}
and
\begin{displaymath}
 d\hat{y}=dy_1 \wedge \cdots \wedge \widehat{dy_b} \wedge\cdots \wedge dy_n.
\end{displaymath}
Integrating 
\begin{displaymath}
 \int_0^{\alpha_b}\frac{\partial\psi}{\partial y_a}(\xi(x,y))G^{ab}\frac{\partial e^{2N\varphi(x,y)}}{\partial y_b}dy_b,
\end{displaymath}
by parts we get
\begin{displaymath}
 -\int_0^{\alpha_b}\frac{\partial}{\partial y_b}\left(\frac{\partial\psi}{\partial y_a}(\xi(x,y))G^{ab} \right)e^{2N\varphi(x,y)}dy_b
\end{displaymath} 
and two boundary terms
\begin{displaymath}
\left[\frac{\partial\psi}{\partial y_a}(\xi(x,y))G^{ab} e^{2N\varphi(x,y)}\right]_{y_b=\alpha_b}
\end{displaymath}
and
\begin{displaymath}
-\left[\frac{\partial\psi}{\partial y_a}(\xi(x,y))G^{ab} e^{2N\varphi(x,y)}\right]_{y_b=0}.
\end{displaymath}
The integral (\ref{integral}) must then be equal to
\begin{displaymath}
 -\int_{\triangle}\frac{\partial}{\partial y_b}\left(\frac{\partial\psi}{\partial y_a}(\xi(x,y))G^{ab} \right)e^{2N\varphi(x,y)}dy_b,
\end{displaymath} 
plus two sums of boundary terms
\begin{displaymath}
\int_{y_b=\alpha_b}\frac{\partial\psi}{\partial y_a}(\xi(x,y))G^{ab} e^{2N\varphi(x,y)}
\end{displaymath}
and \begin{displaymath}
-\int_{y_b=0}\frac{\partial\psi}{\partial y_a}(\xi(x,y))G^{ab} e^{2N\varphi(x,y)}.
\end{displaymath}
The first of these boundary terms is zero because $G^{ab}=0$ when $y_b=\alpha_b$. Now consider the second boundary term. The set where $y_b=0$ contains the facet $F_b$. For this facet, we can take $u_b=e_b$. We know $G^{ab}u_b=0$ on the facet $F_b$, this means that
\begin{displaymath}
 G^{ab}(y)=0, \forall y \in F_b,
\end{displaymath}
and the second boundary term is also zero.$\square$

After the integration by parts, equation (\ref{sharp}) becomes 
\begin{displaymath}
 \psi^\sharp(x)=\psi(x)+\frac{1}{2N}\int_{\triangle}\frac{1}{\int_{\triangle}e^{2N\varphi(x,y)} dy}\frac{\partial }{\partial y_b}\left(\frac{\partial\psi}{\partial y_a}(\xi(x,y))G^{ab}\right)e^{2N\varphi(x,y)}dy.
\end{displaymath}
Now we apply this process again. Set
\begin{displaymath}
 \psi_1(x,y)=\frac{\partial }{\partial y_b}\left(\frac{\partial\psi}{\partial y_a}(\xi(x,y))G^{ab}\right).
\end{displaymath}
We write the Taylor expansion for $\psi_1$ in $y$ around $x$
\begin{displaymath}
 \psi_1(x,y)=\psi_1(x,x)+\frac{\partial\psi_1}{\partial y_a}(\xi_1(x,y))(y_a-x_a),
\end{displaymath}
but we are only interested in the first term of this expansion since the second will bring an $O(\frac{1}{N^2})$ term to equation (\ref{sharp}). Now
\begin{displaymath}
 \psi_1(x,x)=\frac{1}{2}\frac{\partial^2 \psi}{\partial y_a \partial y_b}G^{ab}+\frac{\partial\psi}{\partial y_a}\frac{\partial G^{ab}}{\partial y_b},
\end{displaymath}
where we have used the calculation of the derivatives of $\xi$ at points of the form $(x,x)$. Replacing again in equation (\ref{sharp}) we find 
\begin{displaymath}
 \psi^\sharp(x)=\psi(x)+\frac{1}{2N}\left( \frac{1}{2}\frac{\partial^2 \psi}{\partial y_a\partial y_b}G^{ab}+\frac{\partial\psi}{\partial y_a}\frac{\partial G^{ab}}{\partial y_b}\right)+O\left(\frac{1}{N^2}\right).\square
\end{displaymath} 
We use the equation in the above lemma to substitute for $\psi^\sharp$ in (\ref{somapolytopo}). We get the asymptotic behavior we are interested in, up to terms in $O\left(\frac{1}{N^2}\right)$:
\begin{displaymath}
  \left(1+\frac{1}{N}\left(\frac{\partial}{\partial h_1}+\centerdot\centerdot+\frac{\partial}{\partial h_d}\right)\right) \int_{{\triangle_h}}\left( \psi(x)+\frac{1}{2N}\left( \frac{G^{ab}}{2}\frac{\partial^2 \psi}{\partial y_a \partial y_b}+\frac{\partial\psi}{\partial y_a}\frac{\partial G^{ab}}{\partial y_b}\right)\right),
\end{displaymath}
evaluated at $h=0$. But since $\psi$ has compact support in $\triangle$ we must have
\begin{displaymath}
 \frac{\partial}{\partial h_i}\int_{{\triangle_h}}\psi=0
\end{displaymath}
because this derivative is calculated as the limit, as $h_i$ tends to zero, of the expression
\begin{displaymath}
 \frac{\int_{\triangle_{h_i}\setminus \triangle}\psi}{h_i}
\end{displaymath}
and the numerator is zero since for small enough $h_i$, $\psi$ is zero on the set $\triangle_{h_i} \setminus \triangle$ (note that we may have to consider the set $\triangle\setminus \triangle_{h_i}$ instead). The term in $\frac{1}{N}$ in the expansion of $\int_{\triangle}\psi\mu_N$ is therefore
\begin{displaymath}
 \frac{1}{2}\int_{\triangle}\frac{G^{ab}}{2}\frac{\partial^2 \psi}{\partial y_a \partial y_b}+\frac{\partial\psi}{\partial y_a}\frac{\partial G^{ab}}{\partial y_b}dy,
\end{displaymath}
and we can integrate by parts. Since $\psi$ and its derivatives have compact support in $\triangle$ we do not pick boundary terms. We get 
\begin{displaymath}
 -\frac{1}{4}\int_{\triangle}\psi \frac{\partial^2 G^{ab}}{\partial y_a\partial y_b}.
\end{displaymath}
\subsection{The Tian-Yau-Zelditch asymptotic expansion}
The pointwise asymptotics of the function $\sum |s_m|^2$, where $\{s_m\}$ is an orthonormal basis for the space $H^0(X,L^N)$ was studied in \cite{catlin}, \cite{lu}, \cite{t} and \cite{z}. The following theorem holds:
\begin{theorem}[Catlin,Lu,Tian,Zelditch]
  Let $X$ be a K\"ahler manifold whose symplectic form, $\omega$ has integral cohomology class and $L\rightarrow X$ a line bundle with a Hermitian metric coming from the metric on $X$. Consider an orthonormal basis for the space $H^0(X,L^N)$, $\{s_m\}$. Then there is an asymptotic expansion 
 \begin{displaymath}
 \sum |s_m|^2 \sim A_0(\omega) N^n+A_1(\omega)N^{n-1}+\cdots
\end{displaymath}
and 
\begin{displaymath}
 A_0(\omega)=1, \qquad A_1(\omega)=\frac{s(\omega)}{2},
\end{displaymath}
where $s$ is the scalar curvature. More precisely there exist constants $K_{r}$ such that
\begin{displaymath}
 \left|\left|\sum |s_m|^2-\sum_{i=0}^r  A_i(\omega) N^{n-i}\right|\right|_{\mathcal{C}^0(X)}\leq K_rN^{n-r-1}.
\end{displaymath}
\end{theorem}
Using this result on our toric variety $X$, since our symplectic coordinates are well defined and smooth on our dense open subset, it is easy to conclude that for a compactly supported $\psi$
\begin{displaymath}
 \int_{\triangle}\sum |s_m|^2 \psi dy \sim N^n\int_{\triangle} \psi  +N^{n-1}\int_{\triangle} \psi \frac{s}{2} +\cdots.
\end{displaymath}
We note here that $\sum |s_m|^2$ is actually a function of $y$ only. This is because the torus action on $L$ preserves the Hermitian metric and thus leaves $\sum |s_m|^2$ invariant. Comparing this result with the result obtained in the previous subsection we conclude that 
\begin{displaymath}
 \int_{\triangle} \psi \frac{s}{2}=-\frac{1}{4}\int_{\triangle}\psi \frac{\partial^2 G^{ab}}{\partial y_a\partial y_b}.
\end{displaymath}
Thus, since this holds for all compactly supported $\psi$, we must have
\begin{displaymath}
 s=-\frac{1}{2} \frac{\partial^2 G^{ab}}{\partial y_a \partial y_b},
\end{displaymath}
on $X$, which is Abreu's formula for scalar curvature from \cite{m1}.

\section{Concluding remark}
One could in principle use the very explicit method described in the last section to obtain more terms in the asymptotic expansion of the spectral measure. This would allow one to write down formulas for other $A_i$'s appearing in the Catlin, Lu, Tian, Zelditch Theorem and these formulas would be in terms of polytope data only. For example, it is know (see \cite{lu}) that 
\begin{displaymath}
 A_2(\omega)=\frac{1}{3}\triangle s+\frac{1}{24}(|R|^2-4|Ric|^2+3s^2)
\end{displaymath}
 where $\triangle s$ is the Laplacian of the scalar curvature, $R$ is the curvature tensor and $Ric$ is the Ricci curvature. Therefore, the term in $N^{n-2}$ in the measure asymptotics would allow one to write down an expression for the quantity
 \begin{displaymath}
 |R|^2-4|Ric|^2
\end{displaymath}
in terms of polytope data only.

\end{document}